\documentclass[12pt]{article}
\usepackage{tikz}
\usetikzlibrary{calc,through,backgrounds}
\usepackage{multicol}

\usepackage{bm,amssymb,amsthm,amsfonts,amsmath,latexsym,cite,color, psfrag,graphicx,ifpdf,pgfkeys,pgfopts,xcolor,extarrows}

\newtheorem{theo}{Theorem}[section]

\newtheorem{lem} [theo]{Lemma}
\newtheorem{coro}[theo]{Corollary}
\newtheorem{prop}[theo]{Proposition}

\allowdisplaybreaks

\makeatletter \@addtoreset{equation}{section}
\@addtoreset{theo}{}\makeatother

\usepackage{hyperref}

\def\qed{\hfill \rule{4pt}{7pt}}
\def\pf{\noindent {\it Proof.} }
\def\S{  \mathfrak{S}}

\def\Newton{  \mathrm{Newton}}

\def\F{   \mathcal{F}  }

\begin{document}
\begin{center}
{\Large\bf Vertices of Schubitopes}

\vskip 6mm
{\small  Neil J.Y. Fan$^1$ and Peter L. Guo$^2$ }

\vskip 4mm
$^1$Department of Mathematics\\
Sichuan University, Chengdu, Sichuan 610064, P.R. China
\\[3mm]

$^{2}$Center for Combinatorics, LPMC\\
Nankai University,
Tianjin 300071,
P.R. China

\vskip 4mm

$^1$fan@scu.edu.cn, $^2$lguo@nankai.edu.cn
\end{center}

\noindent{A{\scriptsize BSTRACT}.}
Schubitopes were introduced by  Monical, Tokcan and Yong as
a specific family of generalized permutohedra. It was  proven by Fink,  M\'esz\'aros and  St.$\,$Dizier
that Schubitopes  are the Newton polytopes of the dual characters of
flagged Weyl modules. Important
 cases of Schubitopes include the Newton polytopes of
Schubert polynomials and key polynomials.
In this paper, we develop  a   combinatorial rule to generate  the vertices  of  Schubitopes.
As an application, we show that the vertices of  the Newton polytope of a
key polynomial can be generated by permutations in  a lower interval in the Bruhat order,
settling  a conjecture of
Monical, Tokcan and Yong.

\section{Introduction}

The objective of this paper is to investigate the vertices of
Schubitopes  introduced by  Monical, Tokcan and Yong
\cite{Mon} during their study of Newton polytopes
 in algebraic combinatorics.   Schubitopes are a specific  family of generalized permutohedra  extensively
  studied  by Postnikov \cite{Pos}.  It was conjectured by Monical, Tokcan and Yong
\cite{Mon}  and shown by Fink,  M\'esz\'aros and  St.$\,$Dizier \cite{Fin} that the Newton polytopes of Schubert and key polynomials are   Schubitopes. More generally,
Fink,  M\'esz\'aros and  St.$\,$Dizier \cite{Fin} showed that
 Schubitopes are the Newton polytopes  of the dual characters of
flagged Weyl modules.

We provide a  combinatorial algorithm to generate  the vertices  of  Schubitopes.
As an application, we prove that the vertices of
the Newton polytope of a key polynomial
can be generated  by permutations in  a lower interval in the Bruhat order,
thus confirming  a conjecture of
Monical, Tokcan and Yong \cite[Conjecture 3.13]{Mon}.
This also  establishes a connection between
the Newton polytopes of key polynomials associated to permutations and the
Bruhat interval polytopes introduced by Kodama and Williams \cite{Kod}.

Schibitopes are polytopes associated to diagrams in
an $n\times n$ grid.
A diagram $D$ is a collection of  boxes in
  an $n\times n$ grid. We adopt the notation $[n]=\{1,2,\ldots,n\}$.
We also abbreviate  an $n\times n$ grid to $[n]^2$, and use
 $(i,j)$ to denote  the box   in row $i$ and
column $j$. Here  the rows (resp., columns) are labeled $1,2,\ldots,n$ from top to bottom (resp.,
from left to right).
The Schubitope $\mathcal{S}_D$ associated to $D$ can be defined as follows.
For $1\leq j\leq n$ and  a subset $S$ of,
define a string   $\mathrm{word}_{j, S}(D)$ by reading the $j$-th column of the $n\times n$
grid  from top to bottom and recording:\vspace{-.4cm}
\begin{itemize}
\item $($ if $(i,j)\not\in D$ and $i\in S$;\vspace{-.25cm}
\item $)$ if $(i,j) \in D$ and $i\not\in S$;\vspace{-.25cm}
\item $\star$ if $(i,j) \in D$ and $i\in S$.\vspace{-.4cm}
\end{itemize}
Let
\[\theta_{D}^j(S)=
\text{\#\{paired $()$'s in $\mathrm{word}_{j, S}(D)$}\}+
\text{\#\{$\star$'s in $\mathrm{word}_{j, S}(D)$}\},\]
where the pairing is by the standard ``inside-out'' convention.
For example, for the following diagram and $S=\{1,3\}$, the strings  $\mathrm{word}_{j, S}(D)$
along with the corresponding values  $\theta_D^j(S)$ (which are abbreviated as $\theta^j$)
are illustrated below.

\vspace{-.2cm}
\begin{figure}[h]
\begin{center}
\begin{tikzpicture}

\def\rectanglepath{-- +(4mm,0mm) -- +(4mm,4mm) -- +(0mm,4mm) -- cycle}

\draw (0mm,0mm)--(20mm,0mm)--(20mm,20mm)--(0mm,20mm)--(0mm,0mm);

\draw [step=4mm,dotted] (0mm,0mm) grid (20mm,20mm);

\draw[black,thick] (0mm,4mm) \rectanglepath;
\draw[black,thick] (0mm,16mm) \rectanglepath;
\draw[black,thick] (0mm,0mm) \rectanglepath;
\draw[black,thick] (12mm,0mm) \rectanglepath;
\draw[black,thick] (8mm,0mm) \rectanglepath;
\draw[black,thick] (4mm,8mm) \rectanglepath;
\draw[black,thick] (12mm,8mm) \rectanglepath;
\draw[black,thick] (12mm,12mm) \rectanglepath;
\draw[black,thick] (16mm,12mm) \rectanglepath;

\node at (2mm,18mm) {$\star$};\node at (6mm,18mm) {$($};
\node at (10mm,18mm) {$($};\node at (14mm,18mm) {$($};
\node at (18mm,18mm) {$($};

\node at (14mm,14mm) {$)$};\node at (18mm,14mm) {$)$};

\node at (2mm,10mm) {$($};\node at (6mm,10mm) {$\star$};
\node at (10mm,10mm) {$($};\node at (14mm,10mm) {$\star$};
\node at (18mm,10mm) {$($};

\node at (2mm,6mm) {$)$};

\node at (2mm,2mm) {$)$};\node at (10mm,2mm) {$)$};\node at (14mm,2mm) {$)$};

\node at (50mm,18mm) {$\mathrm{word}_{1, S}(D)=\star())$};
\node at (80mm,18mm) {$\theta^1=2$};
\node at (48.7mm,14mm) {$\mathrm{word}_{2, S}(D)=(\star$};
\node at (80mm,14mm) {$\theta^2=1$};
\node at (49.2mm,10mm) {$\mathrm{word}_{3, S}(D)=(()$};
\node at (80mm,10mm) {$\theta^3=1$};
\node at (50.2mm,6mm) {$\mathrm{word}_{4, S}(D)=()\star)$};
\node at (80mm,6mm) {$\theta^4=2$};
\node at (49.2mm,2mm) {$\mathrm{word}_{5, S}(D)=()($};
\node at (80.5mm,2mm) {$\theta^5=1.$};

\end{tikzpicture}
\end{center}
\end{figure}
\vspace{-.5cm}

Set
\[\theta_D(S)=\sum_{j=1}^n\theta_{D}^j(S).\]
The Schubitope $\mathcal{S}_D$ is  defined by
\[\mathcal{S}_D=\left\{(x_1,\ldots, x_n)\in \mathbb{R}^n\colon \sum_{i\in [n]}  x_i=\#D
\ \  \text{and}\ \ \sum_{i\in S}x_i\leq  \theta_D(S) \ \ \text{for $S\subsetneq [n]$}\right\}.\]
By definition, $\mathcal{S}_D$ is a generalized permutohedren parameterized by $\{\theta_D(S)\}$,
see for example Postnikov \cite{Pos}.


In this paper, we  characterize  the vertices of the Schubitopes $\mathcal{S}_D$ in terms of certain fillings
of $D$.
Let $S_n$
denote the set of permutations of $[n]$. Given a permutation
$w=w_1w_2\cdots w_n\in S_n$,
define $\mathcal{F}_{w}(D)$ to be the filling
of $D$ with the entries of  $w$ as follows.
The filling is described based  on an assignment of the entries of
$w$ into each column of $D$, independent of the order of columns.
For the $j$-th  column $D_j$, fill
the integers $w_1,\ldots,w_n$ in turn into the empty
boxes of $D_j$ as below. From $k=1$ to $k=n$, put $w_k$ into the
 first (from top to bottom) empty box whose row index is larger than
 or equal to $w_k$. If  there are no such empty boxes, then $w_k$ does
 not appear in the filling and skip to $w_{k+1}$.
For example, Figure \ref{diagram} illustrates the filling
$\mathcal{F}_{w}(D)$ for $w=315624$.

\begin{figure}[h]
\begin{center}
\begin{tikzpicture}

\def\rectanglepath{-- +(5mm,0mm) -- +(5mm,5mm) -- +(0mm,5mm) -- cycle}

\draw [step=5mm,dotted] (0mm,0mm) grid (30mm,30mm);
\draw [step=5mm] (0mm,0mm) grid (10mm,5mm);
\draw [step=5mm] (15mm,0mm) grid (25mm,5mm);
\draw (0mm,10mm)--(10mm,10mm);\draw (15mm,0mm)--(15mm,5mm);
\draw [step=5mm] (5mm,5mm) grid (10mm,10mm);
\draw [step=5mm] (0mm,10mm) grid (20mm,15mm);
\draw(10mm,10mm)--(20mm,10mm);
\draw [step=5mm] (0mm,15mm) grid (5mm,20mm);
\draw [step=5mm] (0mm,20mm) grid (15mm,25mm);\draw(5mm,20mm)--(15mm,20mm);

\draw (20mm,25mm) \rectanglepath;
\draw (25mm,15mm) \rectanglepath;
\draw (25mm,10mm) \rectanglepath;
\draw (25mm,5mm) \rectanglepath;\draw (20mm,20mm) \rectanglepath;

\node at (2.5mm,22.5mm) {1};\node at (7.5mm,22.5mm) {1};\node at (12.5mm,22.5mm) {1};\node at (22.5mm,27.5mm) {1};\node at (22.5mm,22.5mm) {2};

\node at (2.5mm,17.5mm) {3};\node at (27.5mm,17.5mm) {3};
\node at (27.5mm,7.5mm) {5};\node at (27.5mm,12.5mm) {1};

\node at (2.5mm,12.5mm) {2};\node at (7.5mm,12.5mm) {3};
\node at (12.5mm,12.5mm) {3};\node at (17.5mm,12.5mm) {3};

\node at (7.5mm,7.5mm) {5};

\node at (2.5mm,2.5mm) {5};\node at (7.5mm,2.5mm) {6};\node at (17.5mm,2.5mm) {1};\node at (22.5mm,2.5mm) {3};

\end{tikzpicture}
\end{center}
\vspace{-6mm}
\caption{The filling $\mathcal{F}_{w}(D)$ for $w=315624$.}
\label{diagram}
\end{figure}

\begin{theo}\label{AAA}
Let $D$ be a diagram of $[n]^2$. Then the vertex set of the Schubitope  $\mathcal{S}_D$  is
\[\{x(w)\colon w\in S_n\},\]
where  $x(w)=(x_1,x_2,\ldots, x_n)$ is the vector such that
 $x_{k}$ $(1\le k\le n)$  is the number of appearances of  $k$ in  $\mathcal{F}_{w}(D)$.
\end{theo}

For the running example as displayed in Figure \ref{diagram}, we have $x(w)=(6,2,6,0,3,1)$.

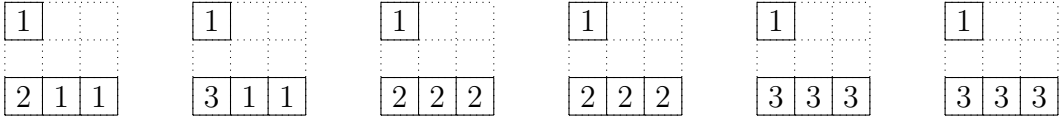
\begin{figure}[h]
\begin{center}
\begin{tikzpicture}

\def\rectanglepath{-- +(5mm,0mm) -- +(5mm,5mm) -- +(0mm,5mm) -- cycle}

\draw [step=5mm,dotted] (0mm,0mm) grid (15mm,15mm);
\draw (0mm,10mm) \rectanglepath;
\draw [step=5mm] (0mm,0mm) grid (15mm,5mm);
\node at (2.5mm,12.5mm) {$1$};
\node at (2.5mm,2.5mm) {$2$};\node at (7.5mm,2.5mm) {$1$};\node at (12.5mm,2.5mm) {$1$};

\draw [step=5mm,dotted] (25mm,0mm) grid (40mm,15mm);
\draw (25mm,10mm) \rectanglepath;
\draw [dotted](25mm,0mm)--(25mm,15mm);
\draw [step=5mm] (25mm,0mm) grid (40mm,5mm);
\draw (25mm,0mm)--(25mm,5mm);
\node at (27.5mm,12.5mm) {$1$};
\node at (27.5mm,2.5mm) {$3$};\node at (32.5mm,2.5mm) {$1$};\node at (37.5mm,2.5mm) {$1$};

\draw [step=5mm,dotted] (50mm,0mm) grid (65mm,15mm);
\draw (50mm,10mm) \rectanglepath;
\draw [dotted](50mm,0mm)--(50mm,15mm);
\draw [step=5mm] (50mm,0mm) grid (65mm,5mm);
\draw (50mm,0mm)--(50mm,5mm);
\node at (52.5mm,12.5mm) {$1$};
\node at (52.5mm,2.5mm) {$2$};\node at (57.5mm,2.5mm) {$2$};\node at (62.5mm,2.5mm) {$2$};

\draw [step=5mm,dotted] (75mm,0mm) grid (90mm,15mm);
\draw (75mm,10mm) \rectanglepath;
\draw [dotted](75mm,0mm)--(75mm,15mm);
\draw [step=5mm] (75mm,0mm) grid (90mm,5mm);
\draw (75mm,0mm)--(75mm,5mm);
\node at (77.5mm,12.5mm) {$1$};
\node at (77.5mm,2.5mm) {$2$};\node at (82.5mm,2.5mm) {$2$};\node at (87.5mm,2.5mm) {$2$};

\draw [step=5mm,dotted] (100mm,0mm) grid (115mm,15mm);
\draw (100mm,10mm) \rectanglepath;
\draw [dotted](100mm,0mm)--(100mm,15mm);
\draw [step=5mm] (100mm,0mm) grid (115mm,5mm);
\draw (100mm,0mm)--(100mm,5mm);
\node at (102.5mm,12.5mm) {$1$};
\node at (102.5mm,2.5mm) {$3$};\node at (107.5mm,2.5mm) {$3$};\node at (112.5mm,2.5mm) {$3$};

\draw [step=5mm,dotted] (125mm,0mm) grid (140mm,15mm);
\draw (125mm,10mm) \rectanglepath;
\draw [dotted](125mm,0mm)--(125mm,15mm);
\draw [step=5mm] (125mm,0mm) grid (140mm,5mm);
\draw (125mm,0mm)--(125mm,5mm);
\node at (127.5mm,12.5mm) {$1$};
\node at (127.5mm,2.5mm) {$3$};\node at (132.5mm,2.5mm) {$3$};\node at (137.5mm,2.5mm) {$3$};

\end{tikzpicture}
\end{center}
\vspace{-6mm}
\caption{The six fillings of $D$ for $D=\{(1,1), (3,1), (3,2), (3,3)\}$.}
\label{6fillings}
\end{figure}

\begin{figure}[h]
\begin{center}
\begin{tikzpicture}
\path [fill=gray] (0mm,0mm)--(40mm,0mm)--(25mm,20mm)--(15mm,20mm)--(0mm,0mm);
\draw[thick](0mm,0mm)--(40mm,0mm)--(25mm,20mm)--(15mm,20mm)--(0mm,0mm);
\node at (-1mm,-3mm) {$\small{(1,0,3)}$};\node at (0mm,0mm) {$\bullet$};
\node at (39mm,-3mm) {$\small{(1,3,0)}$};\node at (40mm,0mm) {$\bullet$};
\node at (28mm,23mm) {$\small{(3,1,0)}$};\node at (25mm,20mm) {$\bullet$};
\node at (12mm,23mm) {$\small{(3,0,1)}$};\node at (15mm,20mm) {$\bullet$};

\node at (7.5mm,10mm) {$\bullet$};
\node at (32.5mm,10mm) {$\bullet$};
\node at (20mm,10mm) {$\bullet$};
\node at (13.33mm,0mm) {$\bullet$};
\node at (26.66mm,0mm) {$\bullet$};

\end{tikzpicture}
\end{center}
\vspace{-6mm}
\caption{The Schubitope $\mathcal{S}_D$ for $D=\{(1,1), (3,1), (3,2), (3,3)\}$.}
\label{nk}
\end{figure}
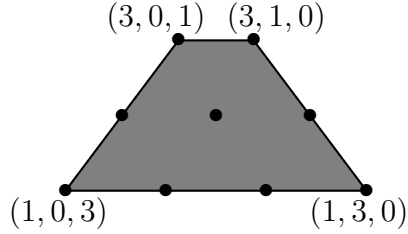

Let us use an example to demonstrate Theorem \ref{AAA}. Let $D=\{(1,1), (3,1), (3,2), (3,3)\}$ be a diagram
of $[3]^2$. The fillings  $\mathcal{F}_{w}(D)$ for the six permutations
$w=123$, $132$, $213$, $231$, $312$, $321$  are listed in Figure \ref{6fillings} in turn
 from left to right.
These fillings generate four vertices:
$x(123)=(3,1,0)$, $x(132)=(3,0,1)$,
$x(213)=x(231)=(1,3,0)$, and $x(312)=x(321)=(1,0,3)$.
The corresponding Schubitope $\mathcal{S}_D$ is a trapezoid as displayed  in Figure \ref{nk},
where the lattice points in $\mathcal{S}_D$ are signified by bullets.

Given a polynomial
$f=\sum_{\alpha\in \mathbb{Z}_{\geq 0}^n}c_\alpha x^\alpha \in \mathbb{R}[x_1,\ldots,x_n],$
  the Newton polytope  of $f$ is the convex hull of the exponent vectors of $f$, namely,
\[\Newton(f)=\mathrm{conv}(\{\alpha\colon c_\alpha\neq 0\}).\]
 Specifying $D$ to the Rothe diagram $D(w)$ of a permutation $w$, the Schubitope  $\mathcal{S}_{D(w)}$
is the Newton polytope  $\Newton(\S_w)$  of the Schubert polynomial $\S_w(x)$
\cite{Mon,Fin}.
Schubert polynomials   were introduced by Lascoux and Sch\"utzenberger \cite{Las}, which
represent the cohomology classes of Schubert cycles in flag varieties.
 Schubert polynomials can be defined in terms of the divided difference operator $\partial_i$,
  which sends
a polynomial  $f $ to
\[\partial_i f =(f
    -s_i f )/ (x_i-x_{i+1}),\]
where $s_i f $ is obtained from $f $ by exchanging $x_i$ and $x_{i+1}$.
For the  permutation $w_0=n \,(n-1)\cdots   1$, set
$\S_{w_0}(x)=x_1^{n-1}x_2^{n-2}\cdots x_{n-1}$.
 For $w\neq w_0$, choose a position $1\leq i<n$ such that $w_i<w_{i+1}$.
Let  $w'$ be the permutation obtained from $w$ by  interchanging
$w_i$ and $w_{i+1}$.
 Set
$\S_{w}(x)=\partial_i \S_{w'}(x)$.

\begin{figure}[h]
\begin{center}
\begin{tikzpicture}

\def\rectanglepath{-- +(5mm,0mm) -- +(5mm,5mm) -- +(0mm,5mm) -- cycle}

\draw [step=5mm,dotted] (50mm,0mm) grid (70mm,20mm);
\draw (50mm,15mm) \rectanglepath;
\draw (50mm,10mm) \rectanglepath;\draw (55mm,10mm) \rectanglepath; \draw (50mm,0mm) \rectanglepath;

\draw[dotted](50mm,0mm)--(50mm,20mm);
\draw [step=5mm,dotted] (0mm,0mm) grid (20mm,20mm);
\draw (5mm,10mm) \rectanglepath;\draw (10mm,10mm) \rectanglepath;
\draw (5mm,5mm) \rectanglepath;
\node at (2.5mm,17.5mm) {$\bullet$};
\node at (17.5mm,12.5mm) {$\bullet$};
\node at (12.5mm,7.5mm) {$\bullet$};
\node at (7.5mm,2.5mm) {$\bullet$};
\draw(2.5mm,17.5mm)--(20mm,17.5mm);\draw(2.5mm,17.5mm)--(2.5mm,0mm);
\draw(17.5mm,12.5mm)--(20mm,12.5mm);\draw(17.5mm,12.5mm)--(17.5mm,0mm);
\draw(12.5mm,7.5mm)--(20mm,7.5mm);\draw(12.5mm,7.5mm)--(12.5mm,0mm);
\draw(7.5mm,2.5mm)--(20mm,2.5mm);\draw(7.5mm,2.5mm)--(7.5mm,0mm);

\node at (10mm,-5mm) {(a)};\node at (60mm,-5mm) {(b)};

\end{tikzpicture}
\end{center}
\vspace{-6mm}
\caption{(a)   $D(w)$ for $w=1432$,   (b)  $D(\alpha)$ for $\alpha=(1,2,0,1)$.}
\label{RS}
\end{figure}
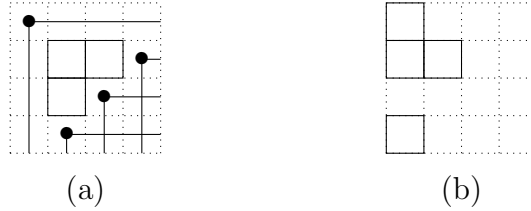

The Rothe diagram
$D(w)$
of   $w\in S_n$  is the diagram
 obtained from the $n\times n$ grid by deleting the box $(i, w_i)$ as well
 as the boxes  to the right of  $(i, w_i)$
or below $(i, w_i)$. Figure \ref{RS}(a)
illustrates the Rothe diagram of $w=1432$.
So, when $D$
 is the Rothe diagram $D(w)$, Theorem \ref{AAA} gives a characterization of the vertices of
  $\Newton(\S_w)$.

When $D$ is restricted to the skyline diagram $D(\alpha)$ of a composition $\alpha$, the Schubitope  $\mathcal{S}_{D(\alpha)}$
is the Newton polytope $\Newton(\kappa_\alpha)$ of the key polynomial $\kappa_\alpha(x)$
\cite{Mon,Fin}.
Key polynomials, also called Demazure characters,
 are characters of the Demazure
modules for the general linear groups \cite{Dem-1,Dem-2}.
Key polynomials  can be defined using the Demazure operator $\pi_i=\partial_ix_i$.
If $\alpha$ is a partition, then set
$\kappa_\alpha(x)=x^\alpha.$
Otherwise, choose $i$ such that $\alpha_i<\alpha_{i+1}$. Let $\alpha'$ be the
composition obtained from $\alpha$ by interchanging $\alpha_i$ and $\alpha_{i+1}$. Set
$\kappa_\alpha(x)=\pi_i \kappa_{\alpha'}(x)$.
It is known that
    $\kappa_\alpha(x)$ can be
 realized as a   specialization of
the nonsymmetric Macdonald polynomial $E_\alpha(x; q, t)$
 at $q=t=\infty$, see Ion \cite{Ion}.
It is also worth mentioning that  every Schubert polynomial is a positive sum
 of key polynomials, see for example Assaf \cite{Ass},
  Lascoux and  Sch\"utzenberger \cite{Las-2}, or Reiner and Shimozono \cite{Rei}.

The skyline
diagram  $D(\alpha)$ of a composition $\alpha$
 is the diagram consisting of the first $\alpha_i$ boxes in row $i$,
 see
Figure \ref{RS}(b)
for the skyline diagram of $\alpha=(1,2,0,1)$.
In this case, Theorem \ref{AAA} can be employed
 to generate the vertices of   $\Newton(\kappa_\alpha)$.

Monical, Tokcan and Yong \cite[Conjecture 3.13]{Mon} conjectured an alternative
characterization  of  the vertices of $\Newton(\kappa_\alpha)$
 in terms of the Bruhat order
on permutations.
Let $\alpha$ be  a composition, and  $\lambda(\alpha)$
be the partition obtained  by resorting the parts of $\alpha$ decreasingly.
Write $w(\alpha)$
for the (unique) permutation of shortest length that sends $\lambda(\alpha)$ to $\alpha$.
Here, given a permutation  $w=w_1\cdots w_n\in S_n$ and a vector $v=(v_1,\ldots, v_n)\in \mathbb{R}^n$,
the (right) action of  $w$ on $v$
is defined as
\begin{equation}\label{ACT}
v\cdot w= (v_{w_1},\ldots, v_{w_n}).
\end{equation}
For two compositions $\alpha,\beta$, define
 \begin{equation}\label{Ord}
 \beta\le \alpha \ \text{ if }\ \lambda(\beta)=\lambda(\alpha) \ \text{ and} \ w(\beta)\le w(\alpha) \ \text{ in the Bruhat order}.
 \end{equation}
 Searles \cite{Sea}  gave an  alternative  description of the partial order in \eqref{Ord}.  For $i<j$ and $\alpha_i<\alpha_j$, let $t_{i,j}(\alpha)$ be
 obtained from $\alpha$ by interchanging $\alpha_i$ and $\alpha_j$.
Then $\beta\leq \alpha$ if and only if $\beta$ can be obtained
 from $\alpha$ by applying a sequence of $t_{i,j}$ \cite[Lemma 3.1]{Sea}.

Based on the decomposition of a key polynomial  into
Demazure atoms \cite{Hag,Las-1,Mas}, Monical, Tokcan and Yong \cite[Theorem 3.12]{Mon}
showed that if $\beta\le\alpha$, then $\beta$ is a vertex of $\Newton(\kappa_{\alpha})$.
They \cite[Conjecture 3.13]{Mon}   conjectured that the converse is still true, that is, if $\beta$ is a vertex of $\Newton(\kappa_{\alpha})$, then $\beta\le\alpha$.
Applying    Theorem \ref{AAA} together with some analysis
on skyline diagrams, we confirm this conjecture.

\begin{theo}\label{conj} Let $\alpha$ be a composition. Then the vertex  set of the Newton polytope   $\Newton(\kappa_{\alpha})$ is  $\{\beta\colon \beta\le\alpha\}$.
\end{theo}

For example, let $\alpha=(1,0,3)$. The key polynomial corresponding to $\alpha$ is
\[\kappa_{\alpha}(x)=x_1^3x_2+x_1^3x_3
+x_1^2x_2^2+x_1^2x_2x_3+x_1^2x_3^2
+x_1x_2^3+x_1x_2^2x_3+x_1x_2x_3^2
+x_1x_3^3.\]
 By Theorem \ref{conj}, it is easily checked that $\Newton(\kappa_{\alpha})$
 has vertex set
 \[\{\beta\colon \beta\le\alpha\}=\{(3, 1, 0), (3, 0, 1), (1, 3, 0), (1, 0, 3)\}.\]
Notice that the  skyline diagram of $\alpha$ is the diagram shown in Figure \ref{6fillings}.
Hence $\Newton(\kappa_{\alpha})$ agrees with the Schubitope in Figure \ref{nk}.

When the parts of $\alpha$ are weakly increasing,   $\kappa_\alpha(x)$
is the Schur polynomial $s_{\lambda(\alpha)}(x)$ \cite{Rei}.
In this case, Theorem \ref{conj}  implies the classical result that the Newton polytope of  a Schur polynomial
 $s_{\lambda}(x)$  is
$\mathcal{P}_{\lambda}$, the permutohedron whose vertices are   rearrangements of $\lambda$.

\noindent
{\it Remark.}
The permutations in $S_n$ are usually
redundant to generate vertices  of a Schubitope, as can be seen in the  example illustrated in Figure \ref{6fillings}.
It is natural to ask which   permutations   are needed to obtain all vertices of
 a Schubitope.
In other words, for two permutations
$w$ and $w'$ in $S_n$, find a characterization to determine whether $x(w)=x(w')$.
Propositions \ref{R-1} and   \ref{R-2} seem  relevant to this question.
When $D$ is a skyline diagram, Theorem \ref{conj} essentially implies that  permutations in a
lower Bruhat interval are enough to generate the vertices.
In the case when $D$ is a Rothe diagram,  we still do not know if Theorem \ref{AAA}
 could be simplified to   a version similar to   Theorem \ref{conj}
 for a skyline diagram.

 Theorem \ref{conj} also  establishes a connection between the Newton polytopes of certain key polynomials  and  Bruhat interval polytopes. For
two permutations $u\leq v$ in the Bruhat order, the    Bruhat interval
polytope $\textsf{Q}_{u,v}$  is the convex hull of the permutations in
the Bruhat interval $[u, v]$.  Bruhat interval polytopes were introduced
by Kodama and Williams \cite{Kod}   in
the context of the Toda lattice and the moment map on the
flag variety, and their combinatorial properties were  studied
by Tsukerman and Williams \cite{Tsu}.
The following corollary is a direct consequence of Theorem \ref{conj}.

\begin{coro}
Let $w=w_1\cdots w_n\in S_n$ be a permutation. View  $w$ as a composition
$(w_1,\ldots, w_n)$. Then the Newton polytope $\Newton(\kappa_{w})$ of $\kappa_w(x)$  is the Bruhat interval polytope $\textsf{Q}_{w, w_0}$, where   $w_0=n \cdots 2 1$
is the largest permutation of $S_n$ in the Bruhat order.
\end{coro}

This paper is structured as follows. In Section 2, we review
a result shown in \cite{Fin} that Schubitopes are Minkowski sums of Schubert matroid polytopes.
This implies that the Schubitope $\mathcal{S}_D$ is a base polytope associated to
some submodular function. Edmonds \cite{Edm} found
 a characterization of vertices of base polytopes for submodular functions.
Based on Edmonds's characterization, we prove Theorem \ref{AAA} in Section 3.
 In the final section, we
 present a proof of Theorem \ref{conj}.

\section{Schubert matroid polytopes}

A matroid is a pair $M=(E, \mathcal{I})$ consisting of a finite set $E$  and a collection
$\mathcal{I}$
 of subsets of $E$, called independent sets, such that
 \begin{itemize}
 \item[(i)] $\emptyset \in \mathcal{I}$;

 \item[(ii)] If $J\in \mathcal{I}$ and $I\subseteq J$, then $I\in \mathcal{I}$;

 \item[(iii)] If $I, J\in \mathcal{I}$ and $|I|<|J|$, then there exists $j\in J\setminus I$
 such that $I\cup \{j\}\in \mathcal{I}$.
 \end{itemize}
 By (ii), a matroid $M$ is determined by the collection $\mathcal{B}$
 of maximal independent sets, called the bases of $M$. So we can   write $M=(E,\mathcal{B})$.
 Moreover, it follows from (iii) that the bases have the same size.
Equivalently, a matroid $M=(E,\mathcal{B})$ can be defined by means of the   exchange  axiom for bases:
\begin{itemize}
 \item[(i')] $\mathcal{B}\neq \emptyset$;

 \item[(ii')] If $A, B\in \mathcal{B} $ and $a\in A\setminus B$, then there exists
 $b\in B\setminus A$ such that $(A\setminus\{a\})\cup \{b\}\in \mathcal{B} $.
 \end{itemize}

Let $S$ be a subset   of $[n]$. The  Schubert matroid $SM_n(S)$ is the matroid with basis
\[\{T\subseteq [n]\colon T\leq S\}.\]
The notation $T\leq S$ means that
\begin{itemize}
\item[(1)]$\# T=\#S$;
\item[(2)] If we write $T=\{a_1<a_2<\cdots<a_k\}$ and
$S=\{b_1<b_2<\cdots<b_k\}$, then $a_i\leq b_i$ for $1\leq i\leq k$.
\end{itemize}
As pointed out by an anonymous referee, Schubert matroids have been rediscovered in different contexts, which
have   been called  freedom matroids, generalized Catalan matroids, PI-matroids,
and   shifted matroids, among others, see Ardila, Fink and   Rinc\'on
 \cite[Example 2.4]{Ard-2}, or the comments
after \cite[Theorem 4.1]{Ard-1} by Ardila
and the comments after \cite[Corollary 3.13]{Bon}  by Bonin,  de Mier and  Noy.
It should also be noted   that  Schubert matroids
are specific families of lattice path matroids
\cite{Bon}, or more generally transversal matroids \cite{Ard-1}
and
positroids \cite[Lemma 23]{Oh}.

 Given a matroid $M=(E, \mathcal{B})$ with $E=[n]$, the associated matroid polytope
 of $M$ is constructed as follows.  Let $\{e_i\colon 1\leq i\leq n\}$ be the   standard basis of $\mathbb{R}^n$.
 For a subset  $B=\{b_1,\ldots, b_k\}$ of $[n]$, write
 \[e_B=e_{b_1}+\cdots+e_{b_k}.\]
 The matroid polytope $P(M)$ is defined by
 \[P(M)=\mathrm{conv}\{e_B\colon B\in \mathcal{B}\}.\]
 The matroid polytope is a generalized permutohedron parameterized   by its rank function
 $\{r_M(S)\}$, see \cite{Fin} for a reference.
 To be specific,
\begin{equation}\label{DE}
P(M)=\left\{x\in \mathbb{R}^n\colon \sum_{i\in [n]}  x_i=r_M([n])
\ \  \text{and}\ \ \sum_{i\in S}x_i\leq  r_M(S)\ \ \text{for $S\subsetneq [n]$}\right\},
\end{equation}
where the rank  function $r_M$  of $M$  is a map from the subsets of $E$ to $\mathbb{Z}_{\geq 0}$
 defined by
 \[r_M(S)=\max\{\# (S\cap B)\colon B\in \mathcal{B}\},\ \ \ \text{for $S\subseteq E$}.\]

It turns out that the
 Schubitope $\mathcal{S}_D$ is the Minkowski sum of Schubert matroid polytopes associated to the columns of $D$. Let  $D$ be a diagram   of $[n]^2$.
 Write $D=(D_1,\ldots, D_n)$,
where, for $1\leq j\leq n$, $D_j$ is the $j$-th column of $D$. The column $D_j$ can be viewed
as a subset of $[n]$:
\[D_j=\{1\leq i\leq n\colon (i,j)\in D\}.\]
Then the column $D_j$ defines a Schubert matroid $SM_n(D_j)$.
For two polytopes $P$ and $Q$, the Minkowski sum of $P$ and $Q$ is defined as
\[P+Q=\{u+v\colon u\in P, v\in Q\}.\]

\begin{theo}[\mdseries{Fink-M\'esz\'aros-St.$\,$Dizier \cite{Fin}}]\label{FFF}
Let $D=(D_1,\ldots, D_n)$ be a diagram of $[n]^2$, and let $r_j$ denote the rank function of
$SM_n(D_j)$. Then
\begin{align}
\mathcal{S}_D&=P(SM_n(D_1))+\cdots+P(SM_n(D_n))\nonumber\\[5pt]
&=\left\{x\in \mathbb{R}^n\colon \sum_{i\in [n]}  x_i=\#D
\ \  \text{and}\ \ \sum_{i\in S}x_i\leq  r_D(S) \ \ \text{for $S\subsetneq [n]$}\right\},
\end{align}
where
\begin{equation}\label{rank}
r_D(S)=r_1(S)+\cdots+r_n(S).
\end{equation}
\end{theo}

\section{Proof of Theorem \ref{AAA}}\label{Sec3}

In this section, we present a proof of Theorem
\ref{AAA}. A crucial observation is that the Schubitope $\mathcal{S}_D$
is the base polytope associated to the function $r_D$. Edmonds \cite{Edm}
obtained a characterization of the vertices of any given base polytope. Based on Edmonds's
characterization, we arrive at a proof of Theorem \ref{AAA}.

\subsection{Schubitopes are base polytopes}

Base polytopes are polytopes associated to submodular functions.
A function $f$ from the subsets of $[n]$ to $\mathbb{R}$ is called a submodular
function, if, for any subsets $S, T\subseteq [n]$,
\[f(S)+f(T)\geq f(S\cup T) + f(S\cap T).\]
To a  submodular function $f$, the  associated base polytope $B_f$ is defined by
\[B_f=\left\{x\in \mathbb{R}^n\colon x_1+\cdots+x_n=f([n]), \ \
\sum_{i\in S}x_i \leq f(S) \ \ \text{for $S\subsetneq [n]$}\right\}.\]

Using the greedy algorithm, Edmonds \cite{Edm}   obtain the following
description of
the vertices of base polytopes for submodular functions,
see also \cite[Theorem 3.22]{Fuj}.

\begin{theo}[\mdseries{\cite{Edm,Fuj}}]\label{Ed}
Let $f\colon 2^{[n]}\rightarrow \mathbb{R}$ be a submodular function.
Then the vertex set of the base polytope $B_f$ is precisely
\[\{x(w)\colon w\in S_n\},\]
where $x(w)=(x_1,\ldots,x_n)$ is the vector in $\mathbb{R}^n$ defined by
\[x_{w_k}=f(\{w_1,\ldots,w_k\})-f(\{w_1,\ldots,w_{k-1}\}).\]
\end{theo}

A fundamental property of a matroid $M$ is  that its rank
function $r_M$   is submodular \cite{Oxl}. Hence the function $r_D$ defined in
\eqref{rank} is submodular. By Theorem \ref{Ed}, we obtain the following
characterization  of the vertex set of a Schubitope.

\begin{theo}\label{EDM}
Let $D$ be a diagram of $[n]^2$.
Then the vertex set of   the Schubitope $\mathcal{S}_D$ is
\[\{x(w)\colon w\in S_n\},\]
where $x(w)=(x_1,\ldots,x_n)$ is the vector in $\mathbb{R}^n$ defined by
\[x_{w_k}=r_D(\{w_1,\ldots,w_k\})-r_D(\{w_1,\ldots,w_{k-1}\}).\]
\end{theo}

\subsection{Rank function  of a  Schubert matroid}\label{sub-1}

Throughout this subsection, we let $C$  be a column of a diagram    of $[n]^2$.
Of course, we can regard $C$ itself as a diagram of $[n]^2$ such that
the boxes lie in exactly one column. Let
$SM_n(C)$ be the Schubert matroid associated to $C$. We show
that the filling $\mathcal{F}_{w}(C)$ generated by the algorithm
 in Introduction can be used to
compute  the rank function $r_C$  of
$SM_n(C)$.
This, together with Theorem \ref{EDM}, leads to a proof of Theorem \ref{AAA}.

A filling $\F$ of $C$ is an assignment of positive integers into some of
the boxes of $C$. A box of $\F$ is called empty if it is not assigned any number.
A filling $\F$ is called column-strict if the numbers appearing in
$\F$ are distinct, and $\F$ is called  flagged if for any nonempty box in row $i$,
the number assigned in it does not
 exceed $i$.
For a subset  $S$ of $[n]$, we denote by $\F(C,S)$  the set of column-strict flagged
fillings $\F$ of $C$ such that all the integers  appearing in $\F$ belong to $S$.
We also denote $\F_{\leq}(C,S)$ to be the subset consisting of the fillings
$\F\in \F(C,S)$ such that the numbers in $\F$ are increasing from top to bottom.
Let $|\mathcal{\F}|$  denote the number of
non-empty boxes of $\F$.

For a permutation $\pi$ of a subset $S$ of $[n]$, we can generate
a filling $\mathcal{F}_{\pi}(C)$ of $C$ by  the algorithm given in Introduction.
Notice that there may exist   empty boxes in $\mathcal{F}_{\pi}(C)$.

\begin{theo}\label{Main-1}
Let $C$ be a column of a diagram of $[n]^2$, and $r_C$ be
the rank function of  $SM_n(C)$. For a $k$-subset $S$ of $[n]$, let $\pi=\pi_1\pi_2\cdots \pi_k$ be any given permutation of elements of $S$.
Then
\begin{equation}\label{ff}
r_C(S)=|\mathcal{F}_{\pi}(C)|.
\end{equation}
\end{theo}

To prove Theorem \ref{Main-1}, we need the following characterization
of the rank function $r_C$.

\begin{theo}\label{T-1}
For any subset $S$ of $[n]$, we have
\begin{equation}
r_C(S)=\max\{|\F|\colon \F\in \F(C,S)\}.
\end{equation}
\end{theo}

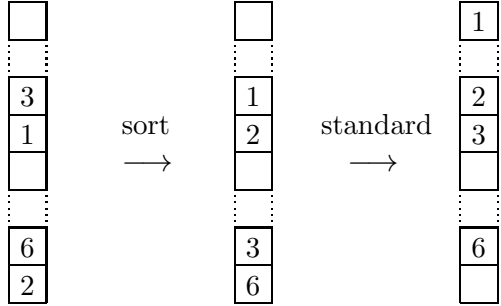
\begin{figure}[h]
\setlength{\unitlength}{0.5mm}
\begin{center}
\begin{picture}(130,80)

\put(0,0){\line(1,0){10}}\put(0,10){\line(1,0){10}}\put(0,20){\line(1,0){10}}
\put(0,0){\line(0,1){20}}\put(10,0){\line(0,1){20}}

\put(0,30){\line(1,0){10}}\put(0,40){\line(1,0){10}}\put(0,50){\line(1,0){10}}\put(0,60){\line(1,0){10}}
\put(0,30){\line(0,1){30}}\put(10,30){\line(0,1){30}}

\qbezier[6](0,20)(0,25)(0,30) \qbezier[6](10,20)(10,25)(10,30)

\qbezier[6](0,60)(0,65)(0,70) \qbezier[6](10,60)(10,65)(10,70)

 \put(0,70){\line(1,0){10}} \put(0,80){\line(1,0){10}}
 \put(0,70){\line(0,1){10}}\put(10,70){\line(0,1){10}}

 \put(3,2.2){\small{2}} \put(3,12.2){\small{6}}
  \put(3,42.2){\small{1}} \put(3,52.2){\small{3}}

\put(60,0){\line(1,0){10}}\put(60,10){\line(1,0){10}}\put(60,20){\line(1,0){10}}
\put(60,0){\line(0,1){20}}\put(70,0){\line(0,1){20}}

\put(60,30){\line(1,0){10}}\put(60,40){\line(1,0){10}}\put(60,50){\line(1,0){10}}\put(60,60){\line(1,0){10}}
\put(60,30){\line(0,1){30}}\put(70,30){\line(0,1){30}}

\qbezier[6](60,20)(60,25)(60,30) \qbezier[6](70,20)(70,25)(70,30)

\qbezier[6](60,60)(60,65)(60,70) \qbezier[6](70,60)(70,65)(70,70)

 \put(60,70){\line(1,0){10}} \put(60,80){\line(1,0){10}}
 \put(60,70){\line(0,1){10}}\put(70,70){\line(0,1){10}}

 \put(63,2.2){\small{6}} \put(63,12.2){\small{3}}
 \put(63,42.2){\small{2}} \put(63,52.2){\small{1}}


 \put(120,0){\line(1,0){10}}\put(120,10){\line(1,0){10}}\put(120,20){\line(1,0){10}}
\put(120,0){\line(0,1){20}}\put(130,0){\line(0,1){20}}

\put(120,30){\line(1,0){10}}\put(120,40){\line(1,0){10}}\put(120,50){\line(1,0){10}}\put(120,60){\line(1,0){10}}
\put(120,30){\line(0,1){30}}\put(130,30){\line(0,1){30}}

\qbezier[6](120,20)(120,25)(120,30) \qbezier[6](130,20)(130,25)(130,30)

\qbezier[6](120,60)(120,65)(120,70) \qbezier[6](130,60)(130,65)(130,70)

 \put(120,70){\line(1,0){10}} \put(120,80){\line(1,0){10}}
 \put(120,70){\line(0,1){10}}\put(130,70){\line(0,1){10}}

\put(123,12.2){\small{6}} \put(123,42.2){\small{3}}
 \put(123,52.2){\small{2}} \put(123,72.2){\small{1}}


 \put(30,35){$\longrightarrow$} \put(30,45){\small{sort}}
 \put(90,35){$\longrightarrow$} \put(83,45){\small{standard}}
\end{picture}
\end{center}
\caption{The sorting operation and standardization operation.}
\label{SS}
\end{figure}

To prove Theorem \ref{T-1}, we define  two operations acting on  $\F(C,S)$ and $\F_{\leq}(C,S)$, respectively.
Let $\F\in \F(C,S)$. The first one is the sorting operation, which transforms $\F$ to a filling $\mathrm{sort}(\F)$ by keeping the empty boxes of $\F$ unchanged and
rearranging the entries  of $\F$   increasingly  from top to bottom.
Figure \ref{SS} gives an example to illustrate the sorting operation.

\begin{prop}\label{B-1}
For $\F\in \F(C,S)$, the filling $\mathrm{sort}(\F)$  belongs to $\F_{\leq}(C,S)$.
\end{prop}

\pf This is trivially  true by treating  a Schubert matroid as
a transversal matroid, see the proof of   \cite[Theorem 2.1]{Ard-1}.
Here, we give a simple verification to make it self-contained.
Obviously, $\mathrm{sort}(\F)$  is column-strict. We   need to verify that
$\mathrm{sort}(\F)$ is flagged. Let $a_1a_2\cdots a_k$ be the word by reading the numbers of $\F$ from top to bottom. Define the inversion number $\mathrm{inv}(\F)$ of
$\F$ to be the number of pairs $(i,j)$ such that $a_i>a_j$.

The proof is by induction on $\mathrm{inv}(\F)$. If $\mathrm{inv}(\F)=0$, then
$\mathrm{sort}(\F)=\F\in \F_\leq(C,S)$. We now consider the case  $\mathrm{inv}(\F)>0$.
Choose $i<j$ such that $a_i>a_j$.
Let $\F'$ be the filling obtained
from $\F$ by interchanging $a_i$ and $a_j$. Clearly, $\mathrm{inv}(F')<\mathrm{inv}(F)$.
We claim that $\F'$ belongs to $\F(C,S)$. This can be seen as follows.
Suppose that $a_i$ lies in the $p$-th row of $\F$, and $a_j$ lies in the $q$-th row of $\F$,
where $p<q$. Since $\F$ is flagged, we have $a_i\leq p$ and $a_j\leq q$.
Combining the facts that $a_i>a_j$ and $p<q$, we reach that
$a_i\leq q$ and $a_j\leq p$. This implies that $\F'$ is   flagged, concluding the claim.
By induction, $\mathrm{sort}(\F')$ belongs to $\F_\leq(C,S)$. Since $\mathrm{sort}(\F)=\mathrm{sort}(\F')$,
we complete the proof.
\qed

The second  operation is  the standardization operation acting on $\F_{\leq}(C,S)$.
Let $\F\in \F_{\leq}(C,S)$. The standardization of $\F$ is the filling $\mathrm{standard}(\F)$    obtained
 by moving upwards the numbers in $\F$  as high as possible subject to the flag condition. More precisely, let $a_1<a_2<\cdots<a_k$ be the integers appearing in $\F$.
 Construct a sequence of fillings  $\F=\F^{(0)}, \F^{(1)},\ldots, \F^{(k)}$
  as follows. For $1\leq t\leq k$, $\F^{(t)}$ is generated from
  $\F^{(t-1)}$ according to the following two cases.
\begin{itemize}
  \item[(1)] The row indices of   empty boxes in $\F^{(t-1)}$ above $a_t$ are all strictly smaller
   than $a_t$. In this case, let $\F^{(t)}=\F^{(t-1)}$;

  \item[(2)] There exist  empty boxes in $\F^{(t-1)}$ above $a_t$ with row indices greater
  than or equal to $a_k$. Let $i_t$ be the smallest such row index. Then $\F^{(t)}$
   is obtained from $\F^{(t-1)} $ by moving $a_t$ up to  the box in row $i_t$.
\end{itemize}
Define $\mathrm{standard}(\F)=\F^{(k)}$. By construction, it is easily seen
 that $\mathrm{standard}(\F)$
belongs to $\F_{\leq}(C,S)$.
See Figure \ref{SS} for an illustration of the standardization operation.

We can now give
a proof of Theorem \ref{T-1}.

\noindent
{\it Proof of Theorem \ref{T-1}.}
Let
\begin{equation}
\overline{r}_C(S)=\max\{|\F|\colon \F\in \F(C,S)\}.
\end{equation}
We first show that $\overline{r}_C(S)\leq r_C(S)$. Suppose that $\F_0\in \F(C,S)$
attains  the maximal cardinality   among all fillings in $\F(C,S)$, namely, $\overline{r}_C(S)=|\F_0|$.   Set
\[\F'_0=\mathrm{standard}(\mathrm{sort}(\F_0)).\]
By Proposition \ref{B-1}, $\F'_0$ belongs to $\F_\leq (C,S)$. Let
$\F''_0$ be the filling of $C$ obtained from $\F'_0$ by assigning each
empty box with its row index. By the construction of the standardization
operator, it is easily checked that $\F''_0$ is a column-strict flagged filling of $C$ such that
the integers  in $\F''_0$ are increasing from top to bottom. Hence the set of integers
in $\F''_0$ forms a base, say  $B_0$, of the Schubert matroid $SM_n(C)$. Moreover,
\[\#(S\cap B_0)\geq |\F_0|,\]
which implies that
\[\overline{r}_C(S)=|\F_0|\leq \#(S\cap B_0)\leq r_C(S).\]

We now verify the reverse direction $\overline{r}_C(S)\geq r_C(S)$. Let $B_0$ be a
base of the Schubert matroid $SM_n(C)$ such that $S\cap B_0$ has the maximal cardinality, that is, $r_C(S)=\#(S\cap B_0)$.
Define a filling $\F_{B_0}$ of $C$ as follows:
Assign the elements of $B_0$
into the boxes of $C$ such that the integers  are increasing from top to
bottom, and then delete the integers  not belonging to $S$.
Since $B_0\leq C$,  $\F_{B_0}$ is a filling in
$\F_\leq(C,S)$. As $|\F_{B_0}|=\#(S\cap B_0)$, we see that
\[\overline{r}_C(S)\geq |\F_{B_0}|= r_C(S).\]
This completes the proof.
\qed

Using Theorem \ref{T-1}, we can finish the proof of Theorem \ref{Main-1}.

\noindent
{\it Proof of Theorem \ref{Main-1}.}  We make induction on the cardinality of $S=\{\pi_1,\ldots,\pi_k\}$.
 Consider the initial
case $k=1$. It is obvious that $r_C(S)=1$ or 0, depending on whether $C$ has a box
  with row index greater than or equal to $\pi_1$. So the equality \eqref{ff} holds.

Assume now that $k\geq 2$ and  \eqref{ff}  is true for $k-1$. Let
\[S'=S\setminus \{\pi_k\}=\{\pi_1,\ldots, \pi_{k-1}\}\]
and
$\pi'=\pi_1\pi_2\cdots \pi_{k-1}.$
Recall that
\[r_C(S')=\max\{\#(S'\cap B)\colon B\in \mathcal{B} \}\ \ \ \ \text{and}\ \ \ \
r_C(S)=\max\{\#(S\cap B)\colon B\in \mathcal{B} \},\]
where $\mathcal{B} $ is the basis of the Schubert matroid $SM_n(C)$.
So we see that
\begin{equation}\label{VV}
r_C(S)=r_C(S')\ \ \ \ \  \text{or}\ \ \ \ \  r_C(S)=r_C(S')+1.
\end{equation}

Keep in mind that $\F_\pi(C)$ is obtained from $\F_{\pi'}(C)$ by putting $\pi_k$ into the topmost empty box of $\F_{\pi'}(C)$ subject to the flag condition.
We conclude the proof by considering  the following cases.

\noindent
Case 1. $\mathcal{F}_{\pi}(C)\neq \mathcal{F}_{\pi'}(C)$. In this case, $|\mathcal{F}_{\pi}(C)|=|\F_{\pi'}(C)|+1$. Since $\F_{\pi}(C)\in \F(C,S)$,
it follows from Theorem \ref{T-1} that
\[r_C(S)\geq |\F_{\pi}(C)|=|\F_{\pi'}(C)|+1.\]
By induction, $r_C(S')=|\F_{\pi'}(C)|$. So $r_C(S)\geq r_C(S')+1$. In view of \eqref{VV}, we have
\[r_C(S)= r_C(S')+1=|\F_{\pi}(C)|,\]
as desired.

\noindent
Case 2. $\F_{\pi}(C)=\F_{\pi'}(C)$.  In this case,
there are no allowable empty boxes in $\F_{\pi'}(C)$ to place $\pi_k$.  There are two subcases.

Case I. There are no empty boxes in $\F_{\pi'}(C)$.
By induction, we have
$r_C(S')=|\F_{\pi'}(C)|=\# C$.
Since $r_C(S)\leq \# C$, it follows from \eqref{VV} that
$r_C(S)=r_C(S')=\# C$,
and hence $r_C(S)=|\F_{\pi}(C)|$.

Case II. There exist empty boxes  in $\F_{\pi'}(C)$, but we cannot put $\pi_k$ into any of
these empty boxes.
Suppose that $l$ is the largest  row index of the empty boxes.
Assume that there are $b$ boxes of $C$ lying strictly below row $l$.
By the construction of $\F_{\pi'}(C)$,
each integer  filled in those $b$ boxes below row $l$ is strictly larger than $l$.
As the box in row $l$ is empty, by the construction of $\F_{\pi}(C)$, we have
$\pi_k>l.$

Assume that $r_C(S')=m$. Let   $\pi_{i_1},\ldots,\pi_{i_m}$
be the elements of $S'$ that are filled in $\F_{\pi'}(C)$.
Again, as the box in row $l$ is empty, by the construction of $\F_{\pi'}(C)$,
it is clear that each integer   in the set
$S'\setminus \{\pi_{i_1}, \pi_{i_2},\ldots,\pi_{i_m}\}$
is strictly  larger than $l$.

We aim to show that $r_C(S)= m$. Suppose to the contrary that $r_C(S)\neq m$.
By \eqref{VV}, we have $r_C(S)= m+1$.
By Theorem \ref{T-1}, there is a filling $\F\in\F(C,S)$
such that $|\F|=m+1$. Notice that $\pi_k$ must belong to $\F$, since otherwise $\F$ is
a filling in $\F(C,S')$ which, together with Theorem \ref{T-1}, would imply that
$r_C(S')\geq m+1$, leading to a contradiction.

Assume that $\pi_{j_1},\ldots,\pi_{j_m}, \pi_k\in S$ are
the integers filled in $\F$.
Notice that  each  integer  in the set
\[\{\pi_{j_1},\ldots,\pi_{j_m}\}\setminus \{\pi_{i_1},\ldots,\pi_{i_m}\}\]
is strictly larger than $l$.
Recall that the integers  filled in those $b$ boxes of $\F_{\pi'}(C)$
 below row $l$ are strictly larger than $l$.
So  $\{\pi_{i_1},\ldots,\pi_{i_m}\}$ contains
exactly  $b$ integers strictly  larger than $l$. Thus $\{\pi_{j_1},\ldots,\pi_{j_m}\}$ contains at least
$b$ integers  strictly larger than $l$. Combining the fact that $\pi_k>l$, the set $\{\pi_{j_1},\ldots,\pi_{j_m},\pi_k\}$
contains at least
$b+1$ integers strictly larger than $l$. However, there are exactly $b$ boxes of $C$ with row indices
strictly larger than $l$. This means  the $m+1$ integers $\pi_{j_1},\ldots,\pi_{j_m},\pi_k$
cannot be filled  into the boxes of $C$ to form a flagged filling, leading to a contradiction.
Thus the assumption that $r_C(S)=m+1$ is false. So we have $r_C(S)=m=|\F_{\pi}(C)|$.
This finishes  the proof.
\qed

\subsection{Proof of Theorem \ref{AAA}}

Using Theorem \ref{EDM} and Theorem \ref{Main-1}, we can now present a proof of Theorem \ref{AAA},
which we restate  below.

\noindent
{\bf Theorem 1.1.}
{\it Let $D$ be a diagram of $[n]^2$. Then the vertex set of the Schubitope  $\mathcal{S}_D$  is
\[\{x(w)\colon w\in S_n\},\]
where  $x(w)=(x_1,x_2,\ldots, x_n)$ is the vector such that
 $x_{k}$ $(1\le k\le n)$  is the number of appearances of  $k$ in  $\mathcal{F}_{w}(D)$.
}

\pf
By  Theorem \ref{EDM} and Theorem \ref{Main-1}, we find that
\begin{align*}
x_{w_k}&=r_D(\{w_1,\ldots,w_k\})-r_D(\{w_1,\ldots,w_{k-1}\})\nonumber\\[5pt]
&=\sum_{j=1}^n r_j(\{w_1,\ldots,w_k\})- \sum_{j=1}^nr_j(\{w_1,\ldots,w_{k-1}\})\nonumber\\[5pt]
&=\sum_{j=1}^n  |\F_{w_1\cdots w_k}(D_j)|- \sum_{j=1}^n  |\F_{w_1\cdots w_{k-1}}(D_j)|\nonumber\\[5pt]
&=|\F_{w_1\cdots w_k}(D)|-|\F_{w_1\cdots w_{k-1}}(D)|.
\end{align*}
Thus $x_{w_k}$ is equal to the number of appearances of $w_k$ in $\F_{w_1\cdots w_k}(D)$. It is obvious that the numbers of appearances of $w_k$ in $\F_{w_1\cdots w_k}(D)$ and in $\F_{w}(D)$ are the same, and so $x_{w_k}$ is equal to the number of appearances of $w_k$ in $\F_{w}(D)$.
This completes the proof. \qed

\section{Proof of Theorem \ref{conj}}

Let us  begin by reviewing the Bruhat order.
We view a permutation
$w=w_1w_2\cdots w_n\in S_n$ as a bijection on $[n]$, that is, $w$ maps $i$ to $w(i)=w_i$.
 As usual, for $1\leq i\leq n-1$, let $s_i=(i, i+1)$ denote the adjacent transposition.
So $ws_i$ is the permutation obtained from $w$ by interchanging $w_i$ and
$w_{i+1}$, while $s_i w$ is obtained by interchanging the values $i$ and $i+1$.
For example, for $w=2143$, we have $ws_2=2413$ but $s_2w=3142$.

 Each
permutation can be written as a product of adjacent transpositions. The length $\ell(w)$ of
a permutation $w$ is the minimum $k$ such that $w=s_{i_1}s_{i_2}\cdots s_{i_k}$, and in this case,
$s_{i_1}s_{i_2}\cdots s_{i_k}$ is called a reduced expression of $w$.
The  (strong) Bruhat order $\leq$ on   $S_n$    is the
closure of the following covering relation.
For
$w, w'\in S_n$, we say that  $w$ covers $w'$    if there exists a
 transposition  $t_{ij}=(i,j)$  such that $w=w'\,t_{ij}$ and $\ell(w)=\ell(w')+1$.
 The Bruhat order can also be characterized by the Subword Property, see for example \cite{Hum}.

\begin{theo}[\mdseries{Subword Property}]\label{Sub}
Let $s_{i_1}s_{i_2}\cdots s_{i_k}$ be any given reduced expression of
a permutation $w$. Then $w'\leq w$ in the Bruhat order if
and only if there exists a subexpression of $s_{i_1}s_{i_2}\cdots s_{i_k}$
that is a reduced expression of $w'$.
\end{theo}

\subsection{A decomposition of the set $\{\beta\colon \beta\leq \alpha\}$}\label{sub-4-1}

Recall that for a composition $\alpha=(\alpha_1,\ldots,\alpha_n)$, $\lambda(\alpha)$
 is the partition obtained  by resorting the parts of $\alpha$ decreasingly, and $w(\alpha)$
is the shortest length permutation such that
\[\lambda(\alpha)\cdot w(\alpha)=\alpha,\]
where the action of a permutation on a vector is as defined in \eqref{ACT}.
The permutation $w(\alpha)$ can be read off directly from $\alpha$ as follows. Let $t_1$ be the largest part of $\alpha$ appearing  in $\alpha$  at positions $l_1<l_2<\cdots<l_{a_1}$ from left to right. Then
put $1,2,\ldots, a_1$ in increasing order at the positions $l_1,l_2,\ldots,l_{a_1}$.
Let $t_2$ be the second largest part of $\alpha$, and $t_2$ appears in $\alpha$  at positions $l_1'<l_2'<\cdots<l_{a_2}'$. Then
put $a_1+1,a_1+2,\ldots, a_1+a_2$ in increasing order at the positions $l_1',l_2',\ldots,l_{a_2}'$.
Repeat the same process for the third largest part of $\alpha$, etc.
For example, for $\alpha=(2,0,1,3,2,0,1)$, we have
 $w(\alpha)=2641375$. We can also construct  $w(\alpha)$  by a recursive
procedure. If $\alpha$ is a partition, then $w(\alpha)$ is the identity permutation.
Otherwise, choose a position $r$ such that $\alpha_r<\alpha_{r+1}$. Let
$\alpha'=\alpha\cdot s_r$. Then
\[w(\alpha)=w(\alpha')\,s_r.\]
The above recursive construction eventually
leads to a reduced expression of $w(\alpha)$.

Let $V(\alpha)$ denote the set appearing in Theorem \ref{conj}:
\[
V(\alpha)=\{\beta\colon \beta\le\alpha\}.
\]

\begin{lem}\label{QI}
For any composition $\alpha$, we have
\[V(\alpha)=\{\lambda(\alpha)\cdot \sigma \colon \sigma\leq w(\alpha)\}.\]
\end{lem}

\pf
By  definition  \eqref{Ord}, it is clear that
$\{\beta\colon \beta\le\alpha\} \subseteq\{\lambda(\alpha)\cdot \sigma \colon \sigma\leq w(\alpha)\}$.
We next verify  the reverse inclusion.
Assume that $\sigma\leq w(\alpha)$ and  $\beta=\lambda(\alpha)\cdot \sigma$. We aim to show that $\beta\leq \alpha$.
In other words, we need to  verify $w(\beta)\leq w(\alpha)$.

Let us first give a description of $w(\beta)$.
Suppose that $\alpha$ has $m$ distinct parts, and that for $1\leq i\leq m$,
 the number of appearances  of  the  $i$-th largest
part is equal to $a_i$. Set $b_0=0$, and $b_i=a_1+\cdots+a_i$ for $1\leq i\leq m$.
 It is easy to check that $w(\beta)$ can be obtained from
$\sigma$ by rearranging the integers in the interval $[b_i+1, b_{i+1}]$ ($0\leq i\leq m-1$)
  increasingly from left to right.

The above  description of  $w(\beta)$ leads to an equivalent characterization of $w(\beta)$.
 It is well known that $S_n$ is
the Coxeter group of type $A_{n-1}$, where $n=b_m$, with generating set $\{s_1,s_2,\ldots,s_{n-1}\}$.
 Let
\[J=\{s_1,s_2,\ldots,s_{n-1}\}\setminus \{s_{b_1}, s_{b_2},\ldots, s_{b_m}\}.\]
Let $(S_n)_J$ denote the parabolic subgroup of $S_n$ generated by $J$, and
let $(S_n)_J \,\sigma$ be the right coset of $(S_n)_J $ with respect to $\sigma$.
Then $w(\beta)$ is the (unique) minimal coset representative of  $(S_n)_J\, \sigma$,
that is, $\ell(s_j\,w(\beta))>\ell(w(\beta))$
for any $s_j\in J$. Hence there is a unique $\tau\in (S_n)_J$
satisfying  that $\sigma=\tau\, w(\beta)$ and
$\ell(\sigma)=\ell(\tau)+\ell(w(\beta))$ \cite[Chapter 1.10]{Hum}.
This implies that the concatenation of any two reduced expressions of $\tau$ and
$w(\beta)$ is a reduced expression of $\sigma$, which, combined with
 the Subword Property  in Theorem \ref{Sub},
yields that $w(\beta)\leq \sigma$. Since $\sigma\leq w(\alpha)$, we have
$w(\beta)\leq w(\alpha)$. This completes the proof.
\qed

By Lemma \ref{QI}, we obtain the following decomposition of $V(\alpha)$.

\begin{prop}\label{LEM-F}
Let $\alpha=(\alpha_1,\ldots, \alpha_n)$ be a composition. Assume that there exists
$1\leq r\leq n-1$ such that $\alpha_r<\alpha_{r+1}$. Let $\alpha'=\alpha\cdot s_r$.
Then
\begin{equation}\label{LEM}
V(\alpha)=V(\alpha')\cup \{v\cdot s_r\colon v\in V(\alpha')\}.
\end{equation}
\end{prop}

\pf To conclude  \eqref{LEM}, by Lemma \ref{QI} it suffices to show that
\begin{equation*}
\{\sigma\colon \sigma\leq w(\alpha)\}=\{\tau\colon \tau\leq w(\alpha')\}\cup \{\tau s_r\colon \tau\leq w(\alpha')\}.
\end{equation*}
This can be easily deduced  from the Subword Property.
Since $\alpha_r<\alpha_{r+1}$, from the arguments above Lemma \ref{QI} it follows
that  that
  $w(\alpha)=w(\alpha')s_r$ and  $\ell(w(\alpha))= \ell(w(\alpha'))+1$.
Let $s_{i_1}\cdots s_{i_k}$ be a reduced expression of $w(\alpha')$. Then
$s_{i_1}\cdots s_{i_k} s_r$ is a reduced expression of $w(\alpha)$.

We first show that
\begin{equation}\label{QO}
\{\sigma\colon \sigma\leq w(\alpha)\}\subseteq \{\tau\colon \tau\leq w(\alpha')\}\cup \{\tau s_r\colon \tau\leq w(\alpha')\}.
\end{equation}
There are two cases.

\noindent
Case 1. $s_r$ is not a (right) descent of $\sigma$, that is, $\ell(\sigma)=\ell(\sigma s_r)-1$. In
this case, any reduced expression  of $\sigma$ does not end with $s_r$. This means that we can
choose a subexpression from   $s_{i_1}\cdots s_{i_k}$ to from a reduced expression of $\sigma$, which, by the Subword Property, implies  $\sigma\leq w(\alpha')$.

\noindent
Case 2. $s_r$ is   a (right) descent of $\sigma$, that is, $\ell(\sigma)=\ell(\sigma s_r)+1$.
Then $s_r$ is not a (right) descent of $\sigma s_r$.
As $\sigma s_r\leq \sigma\leq w(\alpha)$, it follows from Case 1 that $\sigma s_r\leq w(\alpha')$.
Since $\sigma=(\sigma s_r) s_r$, we have $\sigma\in \{\tau s_r\colon \tau\leq w(\alpha')\}$.
This verifies \eqref{QO}.

The reverse set inclusion can be checked in a similar manner, and thus is omitted.
\qed

\subsection{Properties on vertices of $\Newton(\kappa_{\alpha})$}

In this subsection, we use Theorem \ref{AAA} to give two relationships on the vertices of
$\Newton(\kappa_{\alpha})$, which will be used in the proof of Theorem \ref{conj}.

\begin{prop}\label{R-1}
Let $\alpha=(\alpha_1,\ldots,\alpha_n)$ be a composition. Assume that there exists $1\le r\le n-1$ such that $\alpha_r<\alpha_{r+1}$, and that  $w$ is a permutation in $S_n$ such that
$r$ appears before $r+1$ in $w$. Then
\begin{equation}\label{rev}
x(w)=x(s_rw)\cdot s_r.
\end{equation}
\end{prop}

 \pf Write $x(w)=(x_1,\ldots,x_n)$. By Theorem \ref{AAA},   $x_k$ is the number of appearances of $k$ in $\F_w(D(\alpha))$.  Let $D_j$ be the $j$-th column of $D(\alpha)$, which is here viewed as
 a subset $\{i\colon (i,j)\in D_j\}$ of $[n]$. It suffices to prove the following
 claim.

\noindent
Claim. The numbers of appearances of $r$
   and $r+1$ in $\F_w(D_j)$ and $\F_{s_rw}(D_j)$
   are exchanged, while, for $k\neq r, r+1$, the number  of appearances of $k$
 in $\F_w(D_j)$ is the same as the the number  of appearances of $k$ in $\F_{s_rw}(D_j)$.

 For ease of description, for any filling $\F$,
 we use $i\in \F$ to mean that the integer $i$ appears in $\F$.
 To verify the Claim, since $\alpha_r<\alpha_{r+1}$, we have the following three cases.

\noindent
 Case 1. $r\not\in D_j$ and $r+1\not\in D_j$.
In this case, it is easy to check that  $\F_{s_rw}(D_j)$ is obtained from
 $\F_{w}(D_j)$ by replacing $r$ (if any) by $r+1$, and replacing $r+1$ (if any) by $r$.

  \noindent
 Case 2. $r\not\in D_j$ and $r+1 \in D_j$.
 This case is essentially the same as Case 1.

  \noindent
 Case 3. $r\in D_j$ and $r+1 \in D_j$. This case is divided into
 the following subcases.

\noindent
 Subcase I. $r\not\in \F_w(D_j)$ and $r+1\not\in \F_w(D_j)$.
 It is easy to check that $\F_{w}(D_j)=\F_{s_rw}(D_j)$.

 \noindent
Subcase II. $r\notin \F_{w}(D_j)$ and $r+1\in \F_{w}(D_j)$.
 Since $r$ appears before $r+1$ in $w$, this case is impossible
 to occur.

\noindent
Subcase III. $r\in \F_{w}(D_j)$ and $r+1\notin \F_{w}(D_j)$.
In this case, we still have two situations to consider.
\begin{itemize}
\item[(1)] $r$ is not filled in the box $(r,j)$. In this case, it
easy to check that $\F_{s_rw}(D_j)$ is obtained from $\F_{w}(D_j)$
by replacing $r$ with $r+1$.

\item[(2)] $r$ is  filled in the box $(r,j)$. Since $r+1$ does not appear in $\F_{w}(D_j)$,
the box $(r+1,j)$ is filled with an integer, say $w_i$, which is smaller than $r$.
By the construction of  $\F_{w}(D_j)$, $w_i$ must appear after $r$, but before $r+1$.
Hence, when we construct $\F_{s_rw}(D_j)$, the box $(r+1,j)$ is occupied by $r+1$,
the box $(r,j)$ is occupied by $w_i$, and   each box other than $(r,j)$ and $(r+1,j)$ is filled with
the same integer as $\F_{w}(D_j)$. This implies that $\F_{s_rw}(D_j)$ is obtained from
$\F_{w}(D_j)$ by replacing $r$ with $r+1$, and then exchanging the values $r+1$ and $w_i$.
The above arguments are best understood by an example as given in Figure \ref{pro4},
where $w=324615$, $r=4$ and $s_rw=325614.$
\begin{figure}[h]
\begin{center}
\begin{tikzpicture}

\def\rectanglepath{-- +(5mm,0mm) -- +(5mm,5mm) -- +(0mm,5mm) -- cycle}

\draw [step=5mm,dotted] (0mm,0mm) grid (5mm,30mm);
\draw [step=5mm] (0mm,0mm) grid (5mm,25mm);
\node at (2.5mm,22.5mm) {$2$};\node at (2.5mm,17.5mm) {$3$};
\node at (2.5mm,12.5mm) {$4$};\node at (2.5mm,7.5mm) {$1$};
\node at (2.5mm,2.5mm) {$6$};

\draw [step=5mm,dotted] (50mm,0mm) grid (55mm,30mm);
\draw [step=5mm] (50mm,0mm) grid (55mm,25mm);
\draw (50mm,0mm)--(50mm,25mm);
\draw[dotted] (50mm,25mm)--(50mm,30mm);

\node at (52.5mm,22.5mm) {$2$};\node at (52.5mm,17.5mm) {$3$};
\node at (52.5mm,12.5mm) {$1$};\node at (52.5mm,7.5mm) {$5$};
\node at (52.5mm,2.5mm) {$6$};

\node at (3mm,-5mm) {$\F_{324615}(D_j)$};\node at (53mm,-5mm) {$\F_{325614}(D_j)$};

\end{tikzpicture}
\end{center}
\vspace{-6mm}
\caption{An illustration of the proof of Subcase (III)(2).}
\label{pro4}
\end{figure}

\end{itemize}

\noindent
Subcase IV. $r\in \F_{w}(D_j)$ and $r+1\in\F_{w}(D_j)$.
This case is similar to Subcase III.
\begin{itemize}
\item[(1)] $r$ is not filled in the box $(r,j)$. In this case, it
easy to check that $\F_{s_rw}(D_j)$ is obtained from $\F_{w}(D_j)$
by interchanging  $r$ and $r+1$.

\item[(2)] $r$ is  filled in the box $(r,j)$,  and $r+1$ is  filled in the box $(r+1,j)$.
In this case, it is
easy to check that $\F_{s_rw}(D_j)=\F_{w}(D_j)$.

\item[(3)] $r$ is  filled in the box $(r,j)$, but $r+1$ is  filled in a box below
$(r+1,j)$. Since $r+1$ is  filled in a box below
$(r+1,j)$,
the box $(r+1,j)$ is filled with an integer, say $w_i$, which is smaller than $r$.
By the same arguments as those in  Subcase III(2),
 we see that  $\F_{s_rw}(D_j)$ is obtained from
$\F_{w}(D_j)$ by interchanging $r$ with $r+1$, and then interchanging   $r+1$ and $w_i$.
\end{itemize}

The above analysis  allows us to conclude the Claim, and so the proof is complete.
 \qed


\begin{prop}\label{R-2}
Let $\alpha=(\alpha_1,\ldots,\alpha_n)$ be a composition.
 Assume that there exists $1\le r\le n-1$ such that $\alpha_r<\alpha_{r+1}$,
 and  that $w$ is a permutation in $S_n$ such that
$r$ appears before $r+1$ in $w$.
 Let $\alpha'=\alpha\cdot s_r$,  and let $x'(w)$ denote the
 vertex of $\Newton(\kappa_{\alpha'})$ labeled by $w$. Then
\begin{equation}\label{rev-2}
x(w)=x'(w).
\end{equation}
\end{prop}

\pf Let $D_j$ be the $j$-th column of $D(\alpha)$. Write $D'=D(\alpha\cdot s_r)$,
and let $D_j'$ be the $j$-th column of $D(\alpha')$.
If $D_j=D_j'$, it is clear  that $\F_{w}(D_j)=\F_{w}(D_j')$.
If $D_j\neq D_j'$, since $\alpha_r<\alpha_{r+1}$, we
 must have $(r,j)\notin D_j$, $(r+1,j)\in D_j$ and $(r,j)\in D_j', (r+1,j)\notin D_j'$.
Keeping in mind that $r$ appears before $r+1$ in $w$,
it is readily checked that
$\F_{w}(D_j')$ is obtained from $\F_{w}(D_j)$ by moving the box   $(r+1,j)$,
together with the integer filled in the box,  up to row $r$.
This, along with Theorem \ref{AAA}, completes the proof.
\qed

\subsection{Proof of Theorem \ref{conj}}

Based on  Propositions   \ref{LEM-F}, \ref{R-1} and \ref{R-2}, we can now
provide a  proof of Theorem \ref{conj}.

\noindent
{\bf Theorem 1.2.}
{\it Let $\alpha$ be a composition. Then the vertex  set of the Newton polytope   $\Newton(\kappa_{\alpha})$ is  $\{\beta\colon \beta\le\alpha\}$.
}

\pf
Denote by $U(\alpha)$ the vertex set of    $\Newton(\kappa_\alpha)$. By Theorem \ref{AAA},
\begin{align}
U(\alpha)=\{x(w)\colon w\in S_n\}.
\end{align}
As mentioned in Introduction, Monical, Tokcan and Yong \cite[Theorem 3.12]{Mon} showed that
$V(\alpha)\subseteq U(\alpha).$
We finish the proof of Theorem \ref{conj} by proving
\begin{align}\label{L-1}
U(\alpha)\subseteq V(\alpha).
\end{align}

The proof  is by induction on the ``reverse'' inversion number of $\alpha$:
\[\mathrm{rinv}(\alpha)=\#\{1\leq i< j\leq n\colon \alpha_i<\alpha_j \}.\]
We first verify \eqref{L-1} for  the case $\mathrm{rinv}(\alpha)=0$. In this case, $\alpha$ is a partition.
So $w(\alpha)$ is the identity permutation,
and thus
$V(\alpha)=\{\alpha\}$.
On the other hand, it is easy to see that for any permutation $w\in S_n$, $\F_{w}(D)$ is the filling
with boxes in row $k$ ($1\le k\le n$)  filled with $k$.  By Theorem \ref{AAA}, we have
$U(\alpha)=\{\alpha\}$. This verifies \eqref{L-1} for the case $\mathrm{rinv}(\alpha)=0$.

We next consider the case $\mathrm{rinv}(\alpha)>0$. Assume that $r$ is a row index such that $\alpha_r<\alpha_{r+1}$.
Let $\alpha'=\alpha\cdot s_r$.
It is obvious that $\mathrm{rinv}(\alpha')=\mathrm{rinv}(\alpha)-1$.
Let $S^{<}_n$ denote the subset consisting of the permutations $w$ of $S_n$
such that $r$ appears before $r+1$. Let $S_{n}^{>}$ denote the complement
of $S^{<}_n$, namely,
\[S_{n}^{>}=\left\{s_rw\colon w\in S_{n}^{<}\right\}.\]
Write
\[U^{<}(\alpha)=\{x(w)\colon w\in S^{<}_n\}\ \ \ \text{and}\
 \ \ U^{>}(\alpha)=\{x(w)\colon w\in S^{>}_n\}.\]
Then
 \[U(\alpha)=U^{<}(\alpha)\cup U^{>}(\alpha).\]
By Proposition \ref{R-1}, we have
 \begin{align}\label{c1}
 U^{>}(\alpha)=\{v\cdot s_r\colon v\in U^{<}(\alpha)\}.
 \end{align}
By Proposition \ref{R-2}, we have
 \begin{align}\label{c2}
 U^{<}(\alpha)=U^{<}(\alpha')\subseteq U(\alpha').
 \end{align}
Therefore,
 \begin{align*}
 U(\alpha)&=U^{<}(\alpha)\cup U^{>}(\alpha)\\[5pt]
 &=U^{<}(\alpha)\cup\{v\cdot s_r\colon v\in U^{<}(\alpha)\}\ \ \ \ \ \ \ \ \text{(by \eqref{c1})}\\[5pt]
 &\subseteq U(\alpha')\cup \{v\cdot s_r \colon v\in U(\alpha') \}\ \ \ \ \ \ \ \ \ \ \text{(by \eqref{c2})}\\[5pt]
 &\subseteq V(\alpha')\cup \{v\cdot s_r \colon v\in V(\alpha') \} \ \ \ \ \ \ \ \ \ \ \text{(by induction)}\\[5pt]
 &=V(\alpha), \ \ \ \ \ \ \ \ \ \ \ \ \ \ \ \ \ \ \ \ \ \ \ \ \ \ \ \ \ \ \ \  \ \ \ \ \  \text{(by Proposition \ref{LEM-F})}
 \end{align*}
 which proves \eqref{L-1}, as desired.
\qed

\vspace{.2cm} \noindent{\bf Acknowledgments.}
We are grateful to the anonymous referees for valuable comments and suggestions. This work was supported by  the National Natural  Science Foundation of China  (Grant No. 11971250) and Sichuan Science and Technology Program (Grant No. 2020YJ0006).


\begin{thebibliography}{99}

\bibitem{Ard-1}
F. Ardila,
The Catalan matroid,
 J. Combin. Theory Ser. A 104 (2003),  49--62.

\bibitem{Ard-2}
F. Ardila, A. Fink and F.  Rinc\'on,
Valuations for matroid polytope subdivisions,
Canad. J. Math. 62 (2010),  1228--1245.


\bibitem{Ass}
S. Assaf, A generalization of Edelman--Greene insertion for Schubert polynomials,
arXiv:1903.05802.

\bibitem{Bon}
J. Bonin, A. de Mier and M. Noy,
Lattice path matroids: enumerative aspects and Tutte polynomials,
J. Combin. Theory Ser. A 104 (2003),  63--94.

\bibitem{Dem-1}
M. Demazure,
D\'esingularisation des vari\'et\'es de Schubert g\'en\'eralis\'ees,
Ann. Sci. \'Ecole  Norm. Sup.  7 (1974), 53--88.

\bibitem{Dem-2}
M. Demazure,
Une nouvelle formule des caract\'eres,
Bull. Sci. Math. 98 (1974), 163--172.

\bibitem{Edm}
J. Edmonds, Submodular functions, matroids, and certain polyhedra,
 1970 Combinatorial Structures and their Applications
 (Proc. Calgary Internat. Conf., Calgary, Alta., 1969),
 pp. 69--87, Gordon and Breach, New York.




\bibitem{Fin}
A. Fink, K. M\'esz\'aros  and A. St.$\,$Dizier,
Schubert polynomials as integer point transforms of generalized permutahedra,
Adv. Math.  332 (2018), 465--475.


\bibitem{Fuj}
S. Fujishige,
Submodular Functions and Optimization, Second edition,
Annals of Discrete Mathematics, 58, Elsevier B. V., Amsterdam, 2005.


\bibitem{Hag}
 J. Haglund, K. Luoto, S. Mason  and S. van Willigenburg,
 Refinements of the Littlewood-Richardson rule,
 Trans. Amer. Math. Soc.  363 (2011), 1665--1686.



\bibitem{Hum}
J.E. Humphreys,
Reflection Groups and Coxeter Groups,  Cambridge
Studies in Advanced Mathematics, No. 29, Cambridge Univ. Press,
Cambridge, 1990.

\bibitem{Ion}
B. Ion,  Nonsymmetric Macdonald polynomials and Demazure characters,
 Duke Math. J. 116 (2003),  299--318.


\bibitem{Kod}
Y. Kodama and L. Williams,
The full Kostant-Toda hierarchy on the positive flag variety,
Comm. Math. Phys. 335 (2015), 247--283.

\bibitem{Las-1}
A. Lascoux and M.-P. Sch\"utzenberger, Keys \& standard bases, Invariant Theory and Tableaux (Minneapolis, MN, 1988), 125--144, IMA Vol. Math. Appl., 19, Springer, New York, 1990.

\bibitem{Las}
A. Lascoux and M.-P. Sch\"utzenberger, Polyn$ \hat{\mathrm o}$mes de Schubert, C. R. Acad. Sci.
Paris  S\'er. I Math. 294 (1982), 447--450.



\bibitem{Las-2}
A. Lascoux and M.-P. Sch\"utzenberger, Tableaux and non-commuative Schubert
polynomials, Func. Anal. Appl. 23 (1989), 63--64.

\bibitem{Mas}
S. Mason, An explicit construction of type $A$ Demazure atoms,
J. Algebraic Combin. 29 (2009), 295--313.


\bibitem{Mon}
C. Monical, N. Tokcan and A. Yong,
Newton polytopes in algebraic combinatorics,
Selecta Math. (N.S.) 25 (2019), no. 5, Paper No. 66.


\bibitem{Oh}
S. Oh,
Positroids and Schubert matroids,
 J. Combin. Theory Ser. A 118 (2011),  2426--2435.


\bibitem{Oxl}
J. Oxley,  Matroid Theory, Second edition, Oxford Graduate Texts in Mathematics, 21,
Oxford University Press, Oxford, 2011.

\bibitem{Pos}
A. Postnikov,
Permutohedra, associahedra, and beyond,
 Int. Math. Res. Not. IMRN (6) 2009, 1026--1106.

\bibitem{Rei}
V. Reiner and M. Shimozono,
Key polynomials and a flagged
Littlewood-Richardson rule,
J. Combin. Theory  Ser. A 70 (1995), 107--143.


\bibitem{Sea}
D. Searles,
Polynomial bases: positivity and Schur multiplication,
Trans.   Amer. Math. Soc. 73 (2020), 819--847.

\bibitem{Tsu}
E. Tsukerman and L. Williams,
Bruhat interval polytopes,
Adv. Math. 285 (2015), 766--810.


\end{thebibliography}
\end{document}